\documentclass[final,3p,times]{elsarticle}

\usepackage{amstext,amsmath,amsfonts,amssymb,amsbsy,bm,amsthm}
\usepackage{verbatim}
\usepackage{color}
\usepackage{moreverb}
\usepackage{caption}
\usepackage{multirow}
\usepackage{algorithmic}
\usepackage{mathtools}
\usepackage[lined,boxed,commentsnumbered]{algorithm2e}
\def\NN{\hbox{I\kern-.2em\hbox{N}}}
\def\RR{{\mathop{{\rm I}\kern-.2em{\rm R}}\nolimits}}
\def\bfca{\mbox{\boldmath$\alpha$}}
\def\bfcb{\mbox{\boldmath$\beta$}}

\def\Q{{\bf Q}}

\def\f{{\bf f}}

\def\n{{\bf n}}
\def\x{{\bf x}}
\def\y{{\bf y}}

\def\bfeta{\mbox{\boldmath$\eta$}}
\def\bfsigma{\mbox{\boldmath$\sigma$}}
\def\bmu{\mbox{\boldmath$\mu$}}

\def\hB{\hat B}

\def\B{\color{black}}

\newcommand{\be}{\begin{equation}}
\newcommand{\ee}{\end{equation}}
\newcommand{\ba}{\begin{eqnarray}}
\newcommand{\ea}{\end{eqnarray}}

\begin{document}

\begin{frontmatter}



\title{Efficient assembly based on B-spline tailored \\ quadrature rules for the  IgA-SGBEM}


\author[label1]{A.~Aimi}
\ead{alessandra.aimi@unipr.it}
\address[label1]{Dept. of Mathematical, Physical and Computer Science, University of Parma,
Parco Area delle Scienze, 53/A, Parma, Italy}
\author[label2]{F.~Calabr\`o}
\ead{calabro@unicas.it}
\address[label2]{Dept. of Electronic and Information Engineering, University of Cassino and southern Lazio,
Via G. di Biasio 43, Cassino (FR), Italy}
\author[label1]{M.~Diligenti}
\ead{mauro.diligenti@unipr.it}
\author[label3]{M.~L.~Sampoli\corref{cor1}}
\ead{marialucia.sampoli@unisi.it}
\address[label3]{Dept. of Information Engineering and Mathematics, University of
Siena,
Via Roma 56, Siena, Italy}
\author[label4]{G.~Sangalli}
\ead{giancarlo.sangalli@unipv.it}
\address[label4]{Dept. of Mathematics, University of Pavia,
Via Ferrata 1, Pavia, Italy}
\author[label5]{A.~Sestini}
\ead{alessandra.sestini@unifi.it}
\address[label5]{Dept. of Mathematics and Computer Science, University of
Florence,
Viale Morgagni 67, Firenze, Italy }
\begin{abstract}

This paper deals with the discrete counterpart of 2D elliptic model problems rewritten in terms of Boundary Integral Equations. The study is done within the framework of Isogeometric Analysis based on B-splines.
In such  a context, the problem of constructing appropriate, accurate and efficient quadrature rules for the Symmetric Galerkin Boundary Element Method is here investigated.
The new integration schemes, together with row assembly and sum
factorization, are used to build a more efficient strategy to derive
the final linear system of equations.  Key ingredients are weighted
quadrature rules tailored for B--splines, that are constructed to be exact in the whole test space, also with respect to the singular kernel. Several simulations are presented and discussed, showing accurate evaluation of the involved integrals and outlining the
superiority of the new approach in terms of computational cost and elapsed
time with respect to the standard element-by-element assembly.
\end{abstract}

\end{frontmatter}

\section{Introduction}
Boundary Element Methods (BEMs) are an important strategy for the numerical solution of linear partial differential equations appearing in many relevant problems in science and
engineering applications, see \cite{aliabadi, chen}. Through the fundamental solution associated to the considered differential equation, a large class of both exterior and interior elliptic Boundary Value Problems (BVPs) can be reformulated as linear Boundary Integral Equations (BIEs), reducing by one the dimension of the computational domain for the discretization, when compared with Finite Element Methods (FEMs).
For all these problems BEMs can be adopted, offering substantial computational advantages
over other numerical techniques. However, in order to achieve a general efficient numerical implementation, a number of issues has to be carefully addressed. One of the most important consists in the accurate approximation of weakly singular, Cauchy singular
and even hyper-singular integrals over the boundary. Such integrals occur also when the
Symmetric Galerkin Boundary Element Method (SGBEM) \cite{Bonnet} is applied to some BIEs. Indeed, in this case a linear symmetric system of equations is obtained, whose coefficient matrix has entries defined as possibly singular double integrals on the boundary of the problem's domain. Thus, considering the outlined difficulties, it is not surprising that, since the first appearances of SGBEM \cite{Costabel0, maier, Sirtori,
Wendland1}, great effort has been devoted to the efficient and accurate computation of the related Galerkin linear system, as proved by many papers investigating this aspect  (see \cite{ref1, Aimi1} and references therein).
Such techniques were combined to the SGBEM scheme based on Lagrangian bases, always adopting the standard element-by-element assembly procedure, as customary when the Lagrangian
basis is used. This implies that all the double integrals were split
into a sum of integrals on pairs of elements where a local quadrature rule was adopted.

\noindent The advent of Isogeometric Analysis (IgA), \cite{BBCHS06,LibroHughes,ref2} has brought a renewed interest in BEMs and very recently it has been combined for the first time to the SGBEM scheme \cite{ADSS1, ADSS2, Nguyen16}. Papers dealing with problems in acoustics \cite{Coox}, airfoils potential flows \cite{ManSalHel}, Stokes flows \cite{HelArrDeS}, fluid-structure-interactions \cite{HelKieDeSRea} can be found in the literature, using the IgA Galerkin BEM formulation. BEM formulations have been used also to construct computational domains for Galerkin-IgA \cite{Falini}. In order to reduce complexity and gain efficiency, the use of collocation \cite{TauRodHug} or mixed collocation is also common, see e.g. \cite{Ginnis,Gu,Joneidi, Li, Marussig, natara,Peng, politis, Scott, Simp2012, SimpsonScott}.

Following the IgA approach, in \cite{ADSS1} B-splines have been used to represent both the domain geometry and the approximated solution of the problem at hand, giving a significant reduction in the dimension of the discretization space required to attain a fixed accuracy with respect to the standard Lagrangian basis. In such preliminary
implementation of the IgA-SGBEM scheme, the entries of the coefficient matrix and of the right-hand side of the linear system are evaluated by using the accurate quadrature rules introduced in \cite{ref1} and suited for an element-by-element implementation of the
SGBEM scheme. However, such quadrature rules are not capable to take any computational advantage from the higher regularity of B-splines.
The obvious consequence is that the method cannot benefit further from its new isogeometric formulation. In order to avoid this drawback, in the present paper new quadrature rules tailored on B-splines are introduced. The aim is to compute efficiently as well as accurately all the occurring double, possibly singular, integrals, and this is a key point for a new, efficient assembly strategy. The considered quadrature rules are of two kinds.
When no singularity appears a B--spline weighted quadrature is adopted \cite{CMP,
CST_pre}. This is constructed for each B-spline basis function of the approximation space, taken as weight function, in order to be exact in a suitable B-spline space, usually a refinement of the previous one. In presence of singularities, a new idea is
applied. Following the classic approach based on modified moments \cite{CMP,DS90,gautschi}, and using the recurrence relations for B-splines \cite{dB}, we propose a new procedure so that, also in this case, the final quadrature is exact in the test space or in a related refinement. As done in \cite{ABCMS, CST_pre} for the IgA-FEM case, we assembly the final matrix via row loop and sum-factorization. This
allows us to get a remarkable gain in terms of computational costs, as confirmed by the experiments. The proposed strategy is ready to be used for problems equipped by Dirichlet data. For mixed problems, hyper-singular integrals can appear which require some more analysis, for example using coordinate transformation \cite{BEMbook,Tell}, subtracting the singularities \cite{Farutin,Gori,Guiggiani}  or splitting the contributions in order to recover quantities that can be evaluated in closed forms \cite{Aimi1}.

The paper is organized as follows. In Section \ref{sect:BIE} the boundary integral formulation of two model problems is briefly introduced. In Section \ref{sect:discretization} some preliminaries and notations are given, which are necessary in the isogeometric setting for the representation of the domain geometry and of the discretization space; then the double integrals defining the entries of the coefficient matrix and of the right-hand side of the IgA-SGBEM linear system are defined. In Section \ref{sect:quad}, using a unified  formulation, two families of quadrature rules are introduced to approximate the double integrals introduced before. The families are two, since different rules are needed to deal with regular or singular integrals, both occurring in any BEM formulation. The formulas used to deal with non singular integrals are taken from \cite{CST_pre}, and they are briefly summarized in Subsection \ref{sect:quad_reg}. Novel formulas addressing the double singular integrals on fixed nodes and tailored on B-splines are introduced in Subsection  \ref{subsect:singquad},
where a recursive relation for the computation of the necessary modified moments for B-spline functions is presented. These last quadrature rules are also preliminarily tested on some singular integrals whose exact solution is available.
Section \ref{sec:elements} deals with the efficient computation and assembly of the IgA-SGBEM linear system.
The new B-spline--oriented IgA-SGBEM assembly strategy is tested through some numerical examples presented and discussed in Section  \ref{sect:numer}. Finally, Section \ref{sect:concl} reports the research conclusions.

\section{Boundary integral model problems} \label{sect:BIE}
In this work we focus on 2D interior or exterior Laplace model problems on
planar domains, assuming boundary Cauchy data of Dirichlet type. In particular we deal with two different geometries: bounded simply connected domains $\Omega\subset \mathbf{R^2}$, and unbounded domains external to an open limited arc.\\
In the first case, denoting with $\Gamma$ the boundary of $\Omega$, assumed sufficiently regular, we deal with the boundary value problem
\begin{equation} \label{BVP_esse}
\left\{ \begin{array}{ll}
\Delta u=0& {\rm in}\; \Omega\,,\\
u=u_{D}& {\rm on}\; \Gamma\,,
\end{array} \right.
\end{equation}
where $u_D$ is the given boundary datum.\\
Choosing a direct approach \cite{chen}, the boundary integral reformulation of (\ref{BVP_esse}) starts from the representation formula for the solution $u$, i.e.
\begin{equation*} 
u(\x)=-\frac{1}{2\,\pi}\int_{\Gamma}
\hat K(\x,\y) \, q(\y) \, d\gamma_\y+\frac{1}{2 \pi} \, \int_{\Gamma}
\frac{\partial \hat K}{\partial \n_y}(\x,\y) \, u_{\cal D}(\y)\, d
\gamma_\y ,\quad \x \in \Omega
\end{equation*}
where $\x=(x_1,x_2)$,  $\y=(y_1,y_2)$, the
kernel $- \frac{1}{2 \pi} \hat K(\x,\y)$ represents the fundamental solution of the 2D Laplace operator, i.e.
\be \label{kernel} \hat K(\x\,,\,\y) := ln (r)\,, \quad
\mbox{with } \quad r :=\|\x-\y\|_2\,, \ee
and $q:=\frac{\partial u}{\partial \n}$.
Hence, the solution $u$ can be evaluated in any point of the domain, provided that we know the flux $q$ on $\Gamma$. So, with a limiting process for $\x$ tending to $\Gamma$ and using the boundary datum, we obtain
the BIE:\\
\be \label{seconda_BIE}
-\frac{1}{2\,\pi}\int_{\Gamma}
\hat K(\x,\y) \, q(\y) \, d\gamma_\y=\frac{1}{2} \, u_{\cal D}(\x)
-\frac{1}{2 \pi} \, \int_{\Gamma}
\frac{\partial \hat K}{\partial \n_y}(\x,\y) \, u_{\cal D}(\y)\, d
\gamma_\y ,\quad \x \in \Gamma \ ,
\ee
in the boundary unknown $q$.
The boundary integral problem (\ref{seconda_BIE}) can be finally set in the following weak form \cite{Wendland2}:\\
\emph{given} $u_D \in H^{1/2}(\Gamma)$, \emph{find} $q \in H^{-1/2}(\Gamma)$ \emph{such that}
\be \label{weakpb}
{\cal A}(q,p)={\cal F}(p),\quad \forall p\in   H^{-1/2}(\Gamma)
\ee
\emph{where}
\be \label{bilf}
{\cal A}(q,p):=-\frac{1}{2\,\pi}\int_{\Gamma}p(\x) \,\int_{\Gamma}
\hat K(\x,\y) \, q(\y) \, d\gamma_\y \, d\gamma_\x
\ee
\emph{and}
\begin{equation*}
{\cal F}(p):=\int_{\Gamma} p(\x) \,
\Big[\frac{1}{2} \, u_{\cal D}(\x)
-\frac{1}{2 \pi} \, \int_{\Gamma}
\frac{\partial \hat K}{\partial \n_y}(\x,\y) \, u_{\cal D}(\y)\, d
\gamma_\y \Big] \, d\gamma_\x
\end{equation*}
In the second case, still denoting with $\Gamma$ an open limited arc in the plane, we consider the following problem
\begin{equation} \label{BVP_arco}
\left\{ \begin{array}{ll}
\Delta u=0& {\rm in}\; \RR^2\setminus \Gamma,\\
u=u_{D}& {\rm on}\; \Gamma\, .
\end{array} \right.
\end{equation}
The BVP in (\ref{BVP_arco}) can model the electrostatic problem of finding the electric potential around a condenser, whose two faces are so near one another to be considered as
overlapped, knowing the electric potential only on the condenser, see e.g. \cite{bansi}. This example is a classic case where BEM are preferred to FEM: a problem in an infinite domain with an obstacle or a source lying on a curve \cite{Costabel0,Stephan}. \\
Choosing an indirect approach \cite{chen, Stephan}, the BIE coming from the boundary integral reformulation of (\ref{BVP_arco}) reads:
\be \label{prima_BIE}
-\frac{1}{2\,\pi}\int_{\Gamma}
\hat K(\x,\y) \, \varphi(\y)\, d\gamma_\y=u_{\cal D}(\x),\quad \x \in \Gamma\,,
\ee
where the unknown density function $\varphi$ represents the jump of $q$ across $\Gamma$.\\
Boundary integral problems (\ref{prima_BIE}) can be set in a weak form similar to (\ref{weakpb}), where the bilinear form ${\cal A}(\varphi,\psi)$ is defined as in \eqref{bilf} and the right-hand side simplifies in:\\
\begin{equation*} 
{\cal F}(\psi):= \int_{\Gamma} \psi(\x) \,
 u_{\cal D}(\x)
\, d\gamma_\x\,.
\end{equation*}

The price for the simplification of the indirect approach
is that the discrete solution obtained by solving the linear
system does not approximate directly the missing Cauchy data on
the boundary, but just the density function appearing in the chosen
integral representation formula, see e.g. \cite{chen}.

\section{IgA-SGBEM discretization}\label{sect:discretization}
Now, we assume that in both the above model problems, $\Gamma$ is defined as a  parametric curve  with no self intersections, closed in the first case and open in the other. Thus $\Gamma$ is the image of a regular invertible function $\f : I\subset \RR \rightarrow
\Gamma \subset \RR^2$ such that, setting $\f(s) := (f_1(s), f_2(s))\,,$ every
point $\x=(x_1,x_2)\in \Gamma$ can be seen as the image of just one value
$s\in I\,,$ where $I = [a\,,\,b]$ if $\Gamma$ is open and $I=[a\,,\,b),$
otherwise.  In more detail, we assume that $\Gamma$ is parametrically
represented by a function $\f$ defined as follows,
\begin{equation*} 
 \f(s) \,:=\, \sum_{i=1}^N \Q_i B_{i,d}(s) , \quad
s \in I \,,
 \end{equation*}
 where the $\Q_i, \, i=1,\ldots,N\,,$ are  ordered
control points assigned in the plane and defining the shape of
$\Gamma$, and where  $\{B_{i,d}(\cdot),\, i=1,\ldots,N\}$ is the
B-spline basis which spans a space $S$ of piecewise $d$--degree polynomials with respect to an assigned set  $\Delta$ of distinct breakpoints in $I.$
Note that $N$ and the regularity required in $S$ at each inner breakpoint can be a priori established. Actually  the definition of the B-spline basis needs the preliminary introduction of the extended knot vector $T =
\{t_1,\,\ldots, t_d, t_{d+1},\ldots, t_{N+1},
t_{N+2}\,,\ldots,t_{N+d+1}\},$ where $ t_{d+1} = a\,,\,
t_{N+1} = b\,,$ and  $t_1 \le \cdots \le t_{N+d+1}\,.$
Each knot $t_i,\  i={d+2},\ldots,N,$   is a possibly repeated occurrence of an inner breakpoint, while the first and last $d$ knots in $T$ are auxiliary knots characterizing a specific B--spline basis. For brevity, when it is not strictly necessary, in the following we shall refer to it by omitting $d$.
The regularity at a certain inner breakpoint in $\Delta$ of any function in $S$ is fully specified  by prescribing the breakpoint multiplicity (integer between $1$ and $d+1$)  in $T$. Since we want $\f$ at least with a continuous image, such multiplicities are required to be all $\le d$ and  one of them is set to $d$ only when $\Gamma$ has an
angular point which will coincide with a control point.  Extended knot vectors
with {\it multiple} or {\it periodic} auxiliary knots are the more commonly
adopted strategies for completing the definition of $T$, see e.g. \cite{dB}.
 Note that the in Computer Aided Geometric Design (CAGD) the standard way to
represent closed curves relies on periodic extended knot vectors and with the last $d$ control points given by  an ordered repetition of the first ones, see e.g. \cite{handbook}. This clearly means that in practice the dimension of the considered spline space $S$ is \B $N-d$ in this case.\\
In the IgA context, dealing with an isoparametric approach, the same basis or an its
appropriate refinement is considered to generate the discretization space $\hat S$ where we will search the approximate solution of weak BIEs.
Actually spaces with dimension greater than $N$ (or $N-d$ in the closed case)
can be clearly used for the experiments, without abandoning the IgA paradigm,
since the {\it knot insertion procedure} can always be adopted to represent
the boundary in a higher dimension spline space with a desired mesh spacing or
with a reduced regularity at the breakpoints. We remind that such a procedure implies introducing new breakpoints or increasing the multiplicity of the existing ones, see e.g. \cite{Farin93}. \\
In more detail, assuming for notational simplicity that no refinement is done,
when $\Gamma$ is an open arc, we difine the space $\hat S$ as
$$\hat S := span\{\hB_1,\cdots,\hB_N\} \,,$$
where
\begin{equation}\label{lift}
\hB_i(\x) := B_i(\f^{-1}(\x))\,, \quad
\x \in \Gamma\,, \qquad  i=1,\ldots,N\,.
\end{equation}
When $\Gamma$ is a closed curve, the dimension of $\hat S$ must be equal to $N-d.$ Thus we fix $\hat S$ as
 \begin{equation*} 
\hat S :=  span\{\hB_1, \cdots, \hB_{N-d}\} \,,
 \end{equation*}
where $\hB_i,\, i=d+1,\ldots,N-d$ are defined as in (\ref{lift})
while we have
 \begin{equation*} 
\hB_i(\x) \,:= B_i^{(c)}
(\f^{-1}(\x))\,, \quad  \x \in \Gamma\,, \qquad
i=1,\ldots,d\,,
\end{equation*}
with
$$
B_i^{(c)}(t) := \left\{ \begin{array}{ll}  B_i(t) & \mbox{if } t \in I \cap [t_i,t_{i+d+1}]\,, \cr
B_{N-d+i}(t) & \mbox{if } t \in I \cap [t_{N-d+i},t_{N+i+1}]\,, \cr
0 & \mbox{otherwise\,.}    \end{array} \right.
$$
In order to simplify the notation,  in the sequel, when $\Gamma$ is closed, we shall omit the superscript $(c)$ to denote the first $d$ cyclic basis elements.\\
The algebraic reformulation of the IgA-SGBEM scheme leads to a linear system of equations, whose unknowns represents the coefficient of the BIE approximate solution w.r.t. the chosen basis, see e.g. \cite{BEMbook}. In more detail, denoting such a solution with   $\hat{\alpha}(\x)  := \sum_{j=1}^{N_{DoF}} \alpha_j \hat B_j(\x)\,, \, \mbox{with } N_{DoF} := dim(\hat S)\,,$  the resulting linear system has size $N_{DoF}$ and is referred as
\begin{equation} \label{sistlin}
A \,\bfca = \bfcb\,,
\end{equation}
where $A$ is a symmetric positive definite matrix, $\bfca = (\alpha_1,\cdots,\alpha_{N_{DoF}})^T$ is the unknown vector and $\bfcb \in \RR^{N_{DoF}}$ is the vector defining the right-hand side which depends on the given Cauchy data.\\
{
In particular, the entries of the matrix $A$ in (\ref{sistlin})
are, up to the coefficient $-\frac{1}{2\pi}$, double integral of the type (see \eqref{bilf}):
 \be \label{Gmatrix}
 \int_{\Gamma} \hB_i(\x)
\int_{\Gamma} \hat K(\x,\y)\ \hB_j(\y)\ d \gamma_\y \ d\gamma_\x \,,
\ee
while the right-hand side involves the computation, depending on the problem at hand, of the following integrals
\be \label{Kt}
\int_{\Gamma}u_D(\x) \hB_i(\x)\, d \gamma_\x \,, \quad
\ \int_{\Gamma} \hB_i(\x) \int_{\Gamma}
\frac{\partial \hat K}{\partial \n_y}(\x,\y) \, u_D(\y)\, d
\gamma_\y d \gamma_\x\,. \ee
Introducing the scalar variables $s$
and $t$ defined as:
$$ 
s := \f^{-1}(\x)\ , \quad t := \f^{-1}(\y)\, $$
and
the {\it parametric speed} associated with $\Gamma,$
\be \label{J}
J(\cdot) := \sqrt{(f'_1(\cdot))^2+(f'_2(\cdot))^2}\,,
 \ee
in following subsections we rewrite integrals in \eqref{Gmatrix}, \eqref{Kt} over the parametric interval $I$.
}

\subsection{Matrix entries}
{ Referring to the double integral in \eqref{Gmatrix}, we can express it as
 \be \label{final_1} I_K^{(i,j)}:=\int_{I} B_i(s)\, J(s)
\int_{I}  K(s,t)\ B_j(t)\, J(t)  \ dt \ ds\,,
\ee
 where $K(s,t) := \hat K(\f(s),\f(t))\,$.
Let us investigate the nature of the involved kernel. Considering
(\ref{kernel}), we can write
 \[
\hat K(\f(s),\f(t)) =\ln(\|\f(s)-\f(t)\|_2^2)^{1/2}=\frac{1}{2} ln \left(
\frac{\|\f(s)-\f(t)\|_2^2} {(s-t)^2} \right) + ln \left|s-t\right|
\]
Setting
\begin{equation} \label{Rst}
R(s,t):= \frac{\|\f(s)-\f(t)\|_2^2} {(s-t)^2} = \left[
\frac{f_1(s)-f_1(t)}{s-t} \right]^2+\left[
\frac{f_2(s)-f_2(t)}{s-t} \right]^2\,,
\end{equation}
we have that
\[
K(s,t) = K_1(s,t)+K_2(s,t)\,,
\]
where
 $$ 
 K_1(s,t)  \,:=\, \frac{1}{2} ln \left(R(s,t) \right)
\ , \qquad K_2(s,t) \,:=\, ln |s-t| \,. $$
Note that $K_1(s,t)$ can be
defined also for $s=t$ extending its definition by continuity,
since
 \begin{equation}  \label{F1_lim}
 \lim_{t \to s} R(s,t)=J^2(s)\,,
\end{equation}
with $J$ defined as in (\ref{J}).\\
 This splitting of $K$ is useful to separate two contributions, the first coming
from the geometry of $\Gamma$  and the last depending just
from its singular nature. Actually, integral \eqref{final_1} can be evaluated as
$$I_K^{(i,j)}=I_{K_1}^{(i,j)} + I_{K_2}^{(i,j)}$$ with
 \begin{equation} \label{IK12}
I_{K_\ell}^{(i,j)} : = \int_{I} B_i(s) J(s) \int_{I}
K_{\ell}(s,t) \, B_j(t) J(t)\,dt \ ds  \ , \qquad \ell=1,2.
 \end{equation}
Different quadrature rules are needed to compute $I_{K_\ell}$, $\ell=1,2$, since for $\ell=2$  weakly
singular integrals can occur because of the logarithmic nature of the kernel $K_2$.
}

\subsection{Right-hand side elements}
{
Referring to the first integral in \eqref{Kt} we can rewrite it as:
 \begin{equation}  \label{final_b1}
b_1^{(i)}:=\int_{I} B_i(s)\, J(s) \, u_D(s) \ ds\,,
\end{equation}
where $u_D(s)=u_D(\f(s))\,$. Hence we have to evaluate a regular integral.\\
The second double integral in \eqref{Kt} can be expressed as
\be \label{final_b2}
b_2^{(i)} :=\int_{I} B_i(s)\, J(s)\, \int_{I}
\overline{K}(s,t) \, u_D(t) \, \, dt \,ds \,,
\ee
 where $u_D(t)=u_D(\f(t))\,$ and\\
$$\overline{K}(s,t) :=
\frac{\partial \hat K}{\partial \n_\y}(\f(s),\f(t))\,J(t)\,.$$
We note that, when $\f \in C^2(I),$ no singularity occurs in (\ref{final_b2}). Indeed, with some computation, we can  write
$$ 
 \overline{K}(s,t) =\,
\frac{(f_2(t)-f_2(s))\ f'_1(t)-(f_1(t)-f_1(s))\
f'_2(t)}{\|\f(s)-\f(t)\|_2^2} \,,$$
that is,
$$
\overline{K}(s,t)=\frac {1}{R(s,t)}\ \frac{(f_2(t)-f_2(s))\
f'_1(t)-(f_1(t)-f_1(s))\ f'_2(t)}{(t-s)^2}\,, $$
with $R(s,t)$ defined in $(\ref{Rst})$. Then, taking into account (\ref{F1_lim}),
we have
$$
\lim_{t \to s} \overline{K}(s,t)=\frac{1}{2} \frac{f'_2(s)f''_1(s)-f'_1(s)f''_2(s)}{J^2(s)}\,.
$$
\\
On the other hand, without assuming  $\f \in C^2(I)\,, $ the kernel $\bar{K}$ becomes strongly
singular on boundaries with corners and weakly singular on Lyapunov curves, see \cite
[Section 7]{Atkinson}.
}

\section{Novel quadrature rules} \label{sect:quad}
As pointed out in the introduction, in the present paper we explore the construction of the  algebraic counterpart of the IgA-SGBEM scheme, by using two different weighted quadrature strategies for  regular and singular integrals. This means that the usual conditions considered for the construction of the quadrature - the exactness requirements - are imposed with respect to a suitable chosen {\it weight} function.

The developed  formulas have in common the vector of nodes - denoted by $\bfeta$ - but the weights will change according to the exactness requirements,
consisting in  imposing their exactness in a suitable spline space with the respect to the selected weight function.
 Several quadrature rules are necessary, all determined by solving a linear system whose coefficient matrix is in any case a B-splines collocation matrix.
The existence of the rules is ensured by the non singularity
of such a matrix, feature achieved by choosing an $\bfeta$ vector fulfilling the Schoenberg-Whitney conditions, see for example \cite{dB}.
\\
The selected exactness spline space is a possible refinement of the test space $S$,  obtained by uniform subdivision of each  element of $\Delta$ into $N_{ref} \ge 1$ elements, where $N_{ref}>1$ is adopted to improve the quadrature accuracy. This - possibly new - partition of the interval $I$ is denoted by $\bar{\Delta}$. The corresponding spline space can be generated by B-splines, denoted in the following by $\bar{B}_j$, $j=1,\dots,N_E$, with $N_E\approx N_{ref} N_{DoF}$.
\\
Our choice of $\bfeta$, that fulfills the
Schoenberg-Whitney conditions with respect to the extended knot vector associated with $\bar B_j, j=1,\ldots,N_E,$ is the following: in the first and last
element of $\bar{\Delta}$, we take $d+2$ uniformly spaced points; on each inner element we take
midpoints and breakpoints. Hence $\bfeta \in \RR^{N_{quad}}$ with $N_{quad}=2d+2N_{ref}N_h-1,$ where $N_h$ denotes the number of elements of $\Delta$ ($N_h \approx N_{DoF}$ when no multiple inner knot is included in the extended knot vector $T$). Note that this simple choice was used also in \cite{CST_pre} and for BEM with splines in \cite{Arnold}. \\
We detail the construction of the quadrature rules in the following two subsections, respectively related to regular and singular integrals. \B

\subsection{B-spline weighted quadrature rule}\label{sect:quad_reg}

When dealing with the external integrals in (\ref{IK12}) and (\ref{final_b2}),  with the integrals (\ref{final_b1}), or even with the inner integral in (\ref{IK12}) for the kernel $K_1$, regular integrands occur, always with a B-spline factor, denoted in the sequel as $B_i(s)$. In this case, following \cite{CST_pre}, we consider this term as a weight and thus we obtain the quadrature:
\be \label{reg_quad}
{\cal Q}^{(i)} [f] := \sum_{n=1}^{N_{quad}} w^{(i)}_n f(\eta_n) \approx  \int_I f(s) {B}_i (s) \,ds .
\ee
Let us detail the construction introducing some notation\footnote{Notice that we define support as the open set where the function is non-zero.}.
\begin{itemize}
\item  $\mathcal{J}^{(i)}:=\{j$ : $1\leq j \leq N_E$, $\textit{supp}(\bar{B}_j) \cap \textit{supp}({B}_i)\ne \emptyset\}$ ; \quad  $N^{(i)}_E :=\# \mathcal{J}^{(i)}$ (number of active exactness functions);
\item $\mathcal{N}^{(i)}:=\{ n$: $1\leq n \leq N_{quad}$, with $\eta_n \in \textit{supp}( B_i) \}$;\qquad $N^{(i)}_{quad}:=\# \mathcal{N}^{(i)}$ (number of active quadrature nodes).
\end{itemize}
The weights $w_n^{(i)}$, $n \in \mathcal{N}^{(i)}$, are determined by imposing the exactness of the formula on all $\bar{B}_j$ with $j \in \mathcal{J}^{(i)}$, the other weights are set to 0. The local exactness requirements then read:
\be \label{reg_exactness}
{\cal Q}^{(i)} [\bar{B}_j] = \sum_{n\in \mathcal{N}^{(i)}} w^{(i)}_n \bar{B}_j (\eta_n) =  \mu^{(i)}_{j} \qquad \forall j\in \mathcal{J}^{(i)} \,,
\ee
where  $\mu^{(i)}_{j}:=\int_I \bar{B}_j (s) {B}_i (s) \,ds$ can be exactly computed since the integrand is a piecewise polynomial function.
Let us note that the considered definition of $\bfeta$ always ensures that $N^{(i)}_E \leq N^{(i)}_{quad} $. Thus,  the conditions in (\ref{reg_exactness}) lead to a possibly underdetermined linear system with maximum rank, thanks to  the Withney-Schoenberg conditions fulfilled by $\bfeta$.
Actually the system is non squared only when the support of $B_i$ includes (but is not limited to) the first or the last element of $\bar\Delta$, where more nodes are taken. In this case we have chosen to solve the system in the minimum Euclidean norm. Note that $N^{(i)}_E$ can be upper bounded independently from $N_{DoF}$, since  $N^{(i)}_E\leq (1+N_{ref})\ (d+1)\,$. Consequently the size of the linear system in (\ref{reg_exactness}) never becomes prohibitive.

This construction has to be repeated for all the basis functions $B_i\,,\ i=1,\dots N_{DoF}$.
A related Pseudo--code is given in Algorithm \ref{algo_quad_reg}.

\begin{algorithm} \LinesNumbered \SetKwInOut{Input}{Input}
\Input{node vector $\bfeta$, B-splines $B_i$, $i=1,\dots, N_{DoF}$, exactness B-spline evaluations $(\mathbb{B})_{j,n}:= \bar{B}_j(\eta_n)$ $j=1,\dots , N_{E} \,,\ n=1,\dots ,N_{quad}$ }
\For{ $i=1,\dots, N_{DoF}$}{
Extract indexes $\mathcal{J}^{(i)}$ of functions $\bar{B}_j$ having support intersecting the support of $B_i$\;
Calculate integrals ${\mu}^{(i)}_j:= \int_I \bar{B}_j(t) {B}_i (t) \, dt \,,\ j\in \mathcal{J}^{(i)} $\;
Extract indices $\mathcal{N}^{(i)}$ of nodes $\eta_n$ belonging to the support of $B_i$\;
Extract the local collocation matrix  ${\mathbb{B}}^{(i)}= (\mathbb{B})_{\mathcal{J}^{(i)},\mathcal{N}^{(i)}} $ \;
Calculate ${\bf w}^{(i)}:=(w_{\nu}^{(i)})_{\nu\in \mathcal{N}^{(i)}}$ (min. Euclidean norm) solution of ${\mathbb{B}}^{(i)} {\bf w}^{(i)} = \bmu^{(i)} $ \, with $\bmu^{(i)}:={\mu}^{(i)}_j$, $j \in \mathcal{J}^{(i)}$;
 }
 \SetKwInOut{Output}{Output}
\Output{vectors ${\bf w}^{(i)}$, $i=1,\dots, N_{DoF}$.}
\caption{Construction of weighted quadrature rules: regular case}
\label{algo_quad_reg}
\end{algorithm}


\subsection{Singular weighted quadrature rule} \label{subsect:singquad}
When instead we deal with the inner integrals in (\ref{IK12}), for $K=K_2$,  weakly singular integrals can occur, as discussed in Section \ref{sect:discretization}, hence a different formula is used.
We choose to isolate the singular term, as done in \cite{Sloan} and introduce a quadrature:
\be \label{sing_quad}
{\cal Q}^{s} [f] := \sum_{n=1}^{N_{quad}} w^{s}_n f(\eta_n) \approx  \int_I f(t) \  ln|t-s| \,dt .
\ee
In order to fix the weight vector $\mathbf{w}^{s}:=(w_n^{s})_{n=1}^{N_{quad}}$, we impose exactness on the B-spline functions $\bar{B}_j$, $j=1,\ldots,N_E$, that is we require the fulfillment of  the following conditions:
\be \label{sing_exactness}
{\cal Q}^{s} [\bar{B}_j] = \sum_{n=1}^{N_{quad}} w^{s}_n \bar{B}_j (\eta_n) =  \mu_j(s) \quad \forall j=1,\dots , N_E,
\ee
with
\be \label{mod_mom_gen}
\mu_j(s) := \int_I \bar{B}_j (t) ln|t-s| \,dt .
\ee
Then,  $\mathbf{w}^{s}$ is the solution of a linear system of size $ N_{E}\times N_{quad}$ with $N_E \le N_{quad}$ whose  coefficient matrix is a B-splines collocation matrix at the quadrature nodes. We remark that this is a -possibly underdetermined- linear system with maximum rank, thanks to  the Withney-Schoenberg conditions. The right hand side $\bmu^{s}:=(\mu_j(s))_{j=1, \dots , N_E}$, of (\ref{sing_exactness}) is calculated by exploiting the recursive formula for B-splines, see below. Recalling that we are dealing with the inner integrals in (\ref{IK12}), for $K=K_2$, and that the formulas introduced in the previous subsection are used for the outer integrals, the parameter $s$  varies among the entries of the vector $\bfeta$.  So different weight vectors ${\bf w}^{\eta_n} \in \RR^{N_{quad}},\, n=1,\ldots,N_{quad}$, are pre-computed.

Let us derive now the recursive approach which can be use to compute the analytic expressions of the \emph{modified moments} introduced in
(\ref{mod_mom_gen}).
Given a B-spline basis of degree $r$, defined with respect to the partition $\bar\Delta$ of $I$, we set
\begin{equation*}  \label{Iq}
 I_q (B_{j,r} , s) := \int_I ln|t-s|\, t^q\, B_{j,r}(t)\, dt\,, \qquad q \in \NN\,.
\end{equation*}
Then we can write
\be \mu_j(s) = I_0(B_{j,d},s)\,,\label{muj}\ee
since $B_{j,d}=\bar B_j$.

 Taking into account the Cox-De Boor recurrence relation of B-splines \cite{dB}, assuming that $t_j$ is the $j$--th entry of the associated extended knot vector,
$$ B_{j,r}(t) \,=\, \frac{t-t_j} {t_{j+r}-t_j} \, B_{j,r-1}(t) \,+\, \frac{ t_{j+r+1} - t} {t_{j+r+1}-t_{j+1}}\, B_{j+1,r-1}(t)\,, $$
we have that a consequent recurrence relation can be obtained for $I_q (B_{j,r}, s)$:
\be
\label{recIq} I_q (B_{j,r} , s) \,=\,   \frac{ I_{q+1} (B_{j,r-1} , s) - t_j\, I_q (B_{j,r-1} , s) } {t_{j+r}-t_j} \,+\,\frac{ t_{j+r+1}\, I_q (B_{j+1,r-1} , s) - I_{q+1} (B_{j+1,r-1} , s) } {t_{j+r+1}-t_{j+1}}\,.
\ee
We remark that, when multiple knots are taken, if $t_{j+r} = t_j\; (t_{j+r+1}=t_{j+1}),$ the first (second) addend in the right-hand side of (\ref{recIq}) must be set to zero.
In order to compute the desired modified moments \eqref{muj} by means of (\ref{recIq}), we need to compute preliminarily $I_k(B_{i,0},s)$, for $k=0,\ldots,d$ and $i=1,\ldots,N_E+d$.
Now, considering that
$$B_{i,0} (t) = \left\{\begin{array}{ll}
1 & \mbox{if } t_i \le t < t_{i+1}\,, \cr
0 & \mbox{otherwise\,,} \cr
\end{array} \right. $$
we obtain that $I_k(B_{i,0},s) =  \int_{t_i}^{t_{i+1}} ln(|t-s|)\, t^k\, dt\,,$ which clearly is a vanishing quantity if $t_i = t_{i+1}.$ In the opposite case,  using the substitution $z = t-s$, we can write
$$ I_k(B_{i,0},s) = \sum_{j=0}^k  \binom {k} {j} \, s^{k-j} \, \int_{t_i-s}^{t_{i+1}-s} \ln |z|\, z^j\, dz\,. $$
Thus, the procedure starts by computing  the quantities:
$$
\mathcal{I}_i^j(s) \,:=  \left\{\begin{array}{ll} \int_{t_i-s}^{\,t_{i+1}-s} \ln |z|\, z^j\, dz\,, & \mbox{if } t_i  <  t_{i+1} \,\wedge\, [t_i\,,t_{i+1}] \subset [a\,,\,b]\,, \cr
0 & \mbox{otherwise\,} \cr \end{array}  \right. \; i=1,\ldots,N_E+d,\quad j=0,\ldots,d\,.
$$
Note that with some basic analytic computation we can derive the following explicit expression of ${\cal I}_i^j(s)$ for the non trivial case:
$$
\mathcal{I}_i^j(s) \,=\,\displaystyle\frac{z^{j+1}}{j+1} \left. \left(\ln |z| -\frac{1}{j+1}\right)\, \right|_{t_i-s}^{t_{i+1}-s} \  .
$$

The described procedure for the computation of the singular modified moments can be seen as a variant of that introduced in \cite{DS90}, where it was defined for Cauchy singular integrals and used for developing quadrature formulas based on cubic spline interpolation.\\
A related Pseudo--code sketching the main steps required for the computation of ${\bf w}^{s}$ is given in Algorithm \ref{algo_quad_sing}. Notice that this procedure can be used -for example in collocation BEM methods- with respect to a generic vector of points $\bfsigma $, while in our case we will have $\bfsigma = \bfeta$.
\begin{algorithm} \LinesNumbered \SetKwInOut{Input}{Input}
\Input{Vector $\bfsigma$ of abscissae in $I,$ quadrature node vector $\bfeta$,  exactness B-spline evaluations $(\mathbb{B})_{j,n}:= \bar{B}_j(\eta_n) $, $j=1,\dots , N_{E} \,,\ n=1,\dots ,N_{quad}$}
Perform $LU$-factorization of matrix $\mathbb{B}$ \;
\For{ $\nu = 1,\dots,\# {\bfsigma} $}{
Compute $\mu_j( \sigma_{\nu} )$, $ j=1,\dots , N_E$\;
Calculate ${\bf w}^{\sigma_{\nu}}$ as (minimum
Euclidean norm) solution of $\mathbb{B} {\bf w}^{\sigma_{\nu}} = \bmu^{\sigma_{\nu}} $ using matrix factorization \;}
\SetKwInOut{Output}{Output}
\Output{vectors ${\bf w}^{\sigma_{\nu}}$, $\nu = 1,\dots,\# {\bfsigma}$}

 \caption{Construction of quadrature rules: singular case}\label{algo_quad_sing}
\end{algorithm}

 Let us remark that we are proposing a new weighted quadrature where the weight is the singular kernel $\ln |t-s|$ for a fixed value of $s$, and exactness is required on B-splines of a fixed degree.
Following the analysis in \cite{Wendland3}, this requirement is the main ingredient for optimal convergence of the BIE approximate solution, because the quadrature error estimate implies consistency of the overall scheme. Then, it is important to have some indication of the quadrature error when applying the quadrature to a generic function. We consider the approximation of the integral
\be
{\mathrm{I}}^{s}[v] :=\int_{-1}^1  \ln |t-s|\, v( t )\, dt\,,
\ee
for functions $v(t)$ such that these values can be evaluated analytically and for $s=\eta_n$, where $\eta_n$ is one of the quadrature points introduced in Section \ref{sect:quad}. Then, the value ${\mathrm{I}}^{s}[v]$ is approximated by ${\cal Q}^{s} [v]$ as in equation \eqref{sing_quad}.
The error is calculated as follows:

\be \label{ERR}
ERR_{d,N_h}[v]:= \dfrac{ \sum_{n=1}^{N_{quad}} \left( {\cal Q}^{\eta_n} [v] - {\mathrm{I}}^{\eta_n}[v] \right)^2 } {\sum_{n=1}^ {N_{quad}} \left( {\cal Q}^{\eta_n} [v] \right)^2 } \ ,
\ee
where $N_{quad}$ is the number of quadrature points of the external quadrature.
The complete analysis of the error of such rules is beyond the scope of the present paper, we only report here some tests that we have run in order to confirm their good behavior. We have done three kind of tests:
\begin{enumerate}
\item First, we have checked that the exactness requirements are fulfilled. This issue is mainly concerned with the condition number of the collocation matrix: we have chosen a-priori fixed quadrature points that could give bad-conditioned matrices. In all our tests we get machine precision in the numerical computation of the modified moments with the quadrature formula (\ref{sing_quad}). Moreover, we have checked that the integral of the monomial (considered in the whole patch) where exactness is due -namely the monomial $t^{d}$- is computed properly. Also in these cases all tests up to degree $d=5$ give exact result up to machine precision.
\item Then, we have considered the computation of moments where exactness is not required: for quadrature rules where exactness is required on B-splines of degree $d$, we have computed approximated integrals of B-splines of degree $d+1$ and of the monomial $t^{d+1}$. Results are reported in Table \ref{tabunica}.
\item At last, we have considered the integral:
\be \label{testdiff}
{\mathrm{I}}^{s}\left[  \dfrac{ \sqrt{1-t^2} }{ t^2 +25 } \right] = \int_{-1}^1 ln|t-s|\,\dfrac{ \sqrt{1-t^2} }{ t^2 +25 }  \, dt = \pi ln(2) + \dfrac{\pi \sqrt{26} }{5} ln\left(\dfrac{\sqrt{25+s^2}}{5+\sqrt{26}} \right)
\ee
taken from reference \cite{Tsal}. This is an interesting benchmark, and in the cited paper the integral is approximated by a procedure that computes an auxiliary integral via  symbolic computations. On such function our aim is to test monotone convergence confirmed by results available in Table \ref{tabunica}.
\end{enumerate} \B

\begin{table} [th]
\begin{center}
\resizebox{1\textwidth}{!}{
\begin{tabular} {|r|r| c c c c c | }
\hline
 \multicolumn{2}{|r|}{$N_h$} & $ 10$ & {$20$}& {$40$}&{$80$}&{$100$}\\
 \hline
$d=2$ & $v(t)=B_{4,3} (t)$ & {$ 3.46\cdot 10^{-7} $} & {$ 1.41\cdot 10^{-7} $}&{$ 5.62\cdot 10^{-8} $}&{$ 2.08\cdot 10^{-8} $} &{$ 1.49\cdot 10^{-8} $} \\
$3$ & $B_{5,4} (t)$& {$ 2.92\cdot 10^{-8} $} & {$ 1.40\cdot 10^{-8} $} & {$ 7.20\cdot 10^{-9} $} & {$ 3.75\cdot 10^{-9}$} & {$ 3.03\cdot 10^{-9} $} \\
$4$ & $B_{6,5} (t)$ & {$ 5.96\cdot 10^{-10} $} & {$ 2.39\cdot 10^{-10} $} & {$ 1.03\cdot 10^{-10} $} & {$4.40\cdot 10^{-11}$} & {$	3.44\cdot 10^{-11}$} \\
$5$ & $B_{7,6} (t)$ & {$ 2.58\cdot 10^{-11}$} & {$	9.64\cdot 10^{-12}$} & {$	3.92\cdot 10^{-12}$} & {$	1.06\cdot 10^{-10}$} & {$	4.67\cdot 10^{-9}$} \\
\hline \hline
$d=2$ & $v(t)=t^3$ &{$ 3.16\cdot 10^{-6} $} & {$ 1.99 \cdot 10^{-7} $}&{$ 1.24 \cdot 10^{-8} $}&{$ 7.74 \cdot 10^{-10} $} &{$ 3.17 \cdot 10^{-10} $} \\
	  & $t^4$ & {$ 5.21\cdot 10^{-5} $} & {$ 3.25 \cdot 10^{-6} $}&{$ 2.04 \cdot 10^{-7} $}&{$ 1.27 \cdot 10^{-8} $} &{$ 5.22 \cdot 10^{-9} $} \\ \hline
$d=3$ & $v(t)=t^4$ & {$ 1.60\cdot 10^{-5} $} & {$ 1.07 \cdot 10^{-6} $}&{$ 6.99 \cdot 10^{-8} $}&{$ 4.47 \cdot 10^{-9} $} &{$ 1.84 \cdot 10^{-9} $} \\
\hline \hline
$d=2$ & $v(t)=\frac{ \sqrt{1-t^2} }{ t^2 +25 }$ & {$ 6.89\cdot 10^{-4} $} & {$ 2.60 \cdot 10^{-4} $}&{$ 9.45 \cdot 10^{-5} $}&{$ 3.37 \cdot 10^{-5} $} &{$ 2.41 \cdot 10^{-5} $} \\
	 $3$& & {$ 4.05\cdot 10^{-4} $} & {$ 1.50 \cdot 10^{-4} $}&{$ 5.30 \cdot 10^{-5} $}&{$ 1.85 \cdot 10^{-5} $} &{$ 1.31 \cdot 10^{-5} $} \\
	 $4$& & {$ 2.92\cdot 10^{-4} $} & {$ 1.06 \cdot 10^{-4} $}&{$ 3.66 \cdot 10^{-5} $}&{$ 1.25 \cdot 10^{-5} $} &{$ 8.80 \cdot 10^{-6} $} \\
	 $5$& & {$ 2.04\cdot 10^{-4} $} & {$ 7.09 \cdot 10^{-5} $}&{$ 2.44 \cdot 10^{-5} $}&{$ 8.23 \cdot 10^{-6} $} &{$ 5.76 \cdot 10^{-6} $} \\
\hline
\end{tabular}}
\end{center}
\caption{ Error $ERR_{d,N_{h}}[v(t)]$ as defined in \eqref{ERR}. Notice that in the case of functions $B_{j,d}(t)$ the integrating function changes both rowwise, due to degree $d$, and columnwise, due to $N_h$. }   \label{tabunica}
\end{table}

\section{Efficient computation and assembly of IgA-SGBEM  linear system} \label{sec:elements}

 Applying quadrature rules introduced in the previous Section, and referring at first to \eqref{IK12}, we can compute
\begin{equation}\label{IK1e2}
\displaystyle I_{K_1}^{(i,j)}\cong \sum_{n_1=1} ^{N^{(i)}_{quad}}w_{n_1}^{(i)}J(\eta_{n_1}^{(i)})
\sum_{n_2=1}^{N^{(j)}_{quad}}w_{n_2}^{(j)}J(\eta_{n_2}^{(j)})K_1(\eta_{n_1}^{(i)}, \eta_{n_2}^{(j)})
\end{equation}
and
\begin{equation}\label{IK1e2_bis}
\displaystyle I_{K_2}^{(i,j)}\cong
\displaystyle\sum_{n_1=1}^{N^{(i)}_{quad}}w_{n_1}^{(i)}J(\eta_{n_1}^{(i)})
\sum_{n_2=1}^{N_{quad}}w_{n_2}^{(\eta_{n_1}^{(i)})}J(\eta_{n_2})B_j(\eta_{n_2})  = \displaystyle\sum_{n_1=1}^{N^{(i)}_{quad}}w_{n_1}^{(i)}J(\eta_{n_1}^{(i)})
\sum_{n_2\in \mathcal{N}^{(j)}} w_{n_2}^{(\eta_{n_1}^{(i)})}J(\eta_{n_2})B_j(\eta_{n_2}) \,,
\end{equation}
where the last equivalence is due to the local support of $B_j\,,$ so that the number of non vanishing addends  does not depend on $N_{quad}$ and therefor neither on $N_{DoF}.$ \\
For the numerical evaluation of $I_{K_2}^{(i,j)}$, we could distinguish the case $|i-j|> d+1$ where the singularity does not occur, so that the regular quadrature could be used as in \eqref{IK1e2}. Nevertheless, in our implementation, we have used quadrature rule \eqref{IK1e2_bis} also for regular integrals in the construction of matrix $I_{K_2}$, because we have checked that this change of rule not only gives no advantage in terms of computational time, but also gives the same accuracy.\\

For the numerical evaluation of \eqref{final_b1}, we can simply consider\\
$$
\displaystyle b_1^{(i)}\cong \sum_{n=1} ^{N^{(i)}_{quad}}w_{n}^{(i)}J(\eta_{n}^{(i)})u_D(\eta_{n}^{(i)})\,.
$$
At last,  since we assume ${\bf f} \in C^2(I),$ for \eqref{final_b2} we proceed by computing
\begin{equation*}
b_2^{(i)}\cong \sum_{n_1=1}^ {N^{(i)}_{quad}}w_{n_1}^{(i)}J(\eta_{n_1}^{(i)})
\sum_{n_2=1}^ {N^{GQ}} w^{GQ}_{n_2} \overline{K}(\eta_{n_1}^{(i)},\eta^{GQ}_{n_2}) u_D(\eta^{GQ}_{n_2}) \,,
\end{equation*}
where $\left(\bfeta^{GQ}, {\bf w}^{GQ} \right)$ is a gaussian quadrature, being the integrand function regular. 

\subsection{Assembly: computational cost} \label{sect:pseudo}

Following what described in the previous sections, the proposed procedure for the assembly of the algebraic counterpart of the IgA-SGBEM scheme can be sketched as follows:
\begin{enumerate}
\item Fix $N_{ref}$ and a choice of $N_{quad}$ quadrature points $\bfeta$.
\item Construct the family of quadrature rules $\left( \eta^{(i)}_n , w^{(i)}_n \right)_{n=1, \dots N^{(i)}_{quad}}, i=1, \dots , N_{DoF} $ following the procedure in Algorithm \ref{algo_quad_reg}. This is done by solving $N_{DoF}$ linear systems of dimension $N_E^{(i)} \le (1+N_{ref})(d+1)$.
\item Construct the family of quadrature rules $\left( \eta_n , w^{\eta_m}_n \right)_{n=1, \dots N_{quad}}, m=1, \dots , N_{quad} $ following the procedure in Algorithm \ref{algo_quad_sing}. This is done by solving $N_{quad}$ linear systems of dimension $N_E \times N_{quad} \approx N_{ref} N_{DoF} \times 2 N_{ref} N_{DoF}$ where only the right  hand size changes, so that the matrix can be factorized just once. We have performed in our tests the usual $LU$ factorization as implemented in Matlab.
\item Assembly the IgA-SGBEM matrix and right-hand side of interest. For computational efficiency, we propose row assembly via sum factorization. Details of such construction can be found in \cite{ABCMS,CST_pre}. In Algorithm \ref{algo_assebly} we present the technique in the case of \eqref{IK1e2}.
\end{enumerate}

\begin{algorithm} \LinesNumbered \SetKwInOut{Input}{Input}
\Input{Quadrature rules}
Set $\mathcal{C}^{(0)}_{n_1,n_2} := K_1(\eta_{n_1},\eta_{n_1}) $, $n_1,n_2=1,\dots , N_{quad} $\;
\For{$i=1,\dots, N_{\mathsf{DoF}}$}{
Compute $\mathcal{C}^{(1)}_{n_1,j} := \sum_{n_2 \in \mathcal{N}^{(j)} }  w_{n_2}^{(j)} J(\eta^{(j)}_{n_2} ) \mathcal{C}^{(0)}_{n_1,n_2}$, $j\in \mathcal{J}^{(i)}$\;
Compute $\mathcal{C}^{(2)}_{i,j} := \sum_{n_1 \in \mathcal{N}^{(i)} }  w_{n_1}^{(i)} J(\eta^{(i)}_{n_1} ) \mathcal{C}^{(1)}_{n_1,j}$\;
Fix $\mathcal{C}^{(2)}_{i,\cdot}$ to obtain the row $I_{K_1}^{(i,\mathcal{J}^{(i)} )}$\;
}
\SetKwInOut{Output}{Output}
\Output{matrix $I_{K_1}.$}
\caption{Construction of matrix $I_{K_1}$ of equation \eqref{IK1e2} by sum-factorization}\label{algo_assebly}
\end{algorithm}
%
Due to the savings introduced by the patch-strategy for the construction of the weighted quadratures, the number of quadrature points  for each parametric variable is $2d+2 N_{ref} N_h -1$. Thus the computational cost for the assembly of the IgA-SGBEM linear system matrix with the introduced B-spline based quadrature strategy and sum factorization is
$O([2d+2 N_{ref} N_h -1]^2)$ function evaluations. 

\bigskip
\noindent {\bf Remark.} In \cite{ADSS1} the  matrix
entries were numerically evaluated within the framework of the
standard element-by-element assembly phase.
{This implies that, within a double nested cycle over the $N_h$ mesh elements of the partition $\Delta$ on the parametrization interval $I$, having indicated by $e^{(i)}_{\ell_i},\, \ell_i=1,\cdots,d+1$ the elements constituting the support of the B-spline $B_i$, we computed and assembled \eqref{final_1} as
 \be \label{elbyelint} I_K^{(i,j)}\cong \sum_{\ell_i=1}^{d+1}\sum_{\ell_j=1}^{d+1}\int_{e^{(i)}_{\ell_i}} B_i(s)\, J(s)
 \int_{e^{(j)}_{\ell_j}}  K(s,t)\ B_j(t)\, J(t)  \ dt \ ds\,,
\ee
taking into account the polynomial nature of the B-splines over each element of their support and using quadrature schemes
introduced in \cite{Aimi1}, suitable for standard local Lagrangian basis functions. In particular, in presence of
kernel singularity, arising for double integration over couples $(e^{(i)}_{\ell_i}, \, e^{(j)}_{\ell_j})$ of coincident or consecutive elements, related double integrals in \eqref{elbyelint}
}
were split into the sum of
a regular part and a singular part: the first was treated using
Gauss-Legendre rule for both inner and outer integrals, while the
other term was evaluated using an interpolatory product rule,
absorbing kernel singularity into the weights, for the inner
integration and the Double Exponential rule (briefly, DE-rule), suitable for integrand function
having weak singularities at the endpoints, for the outer
integration. When the double integration occurred on couples of
"far" elements Gauss-Legendre rule for both inner and outer
integrals was employed. The interested reader is referred to
\cite{Aimi1} and related references for more details on the above resumed quadrature
schemes.\\ Hence, having denoted by
$N_G$ the number of the
Gauss-Legendre rule points, by $N_{DE}$ the number of DE-rule points
and by $N_{prod}$ the number of points of product rules for
singular integration, the computational cost of the element-by-element procedure
for the generation of the IgA-SGBEM linear system matrix, in terms of integrand
function evaluation, was $
(d+1)^2 (N_h^2\,N_G^2+3\,N_h\,N_{DE}\,N_{prod})\,.$


\section{Numerical examples} \label{sect:numer}
{In this section we present and discuss some numerical results obtained using both the old quadrature strategy with the element-by-element assembly and the new B-spline based quadrature with the related new assembly. The old code is developed in Fortran language while the new one is written in Matlab. Simulations have been performed on a laptop equipped by Intel Core i5 CPU (2.53 GHz, 4 Gb RAM, 64 bit OS).}\\
In all the following experiments, the error
$E_M$ represents the maximum of the discrete error on a uniform mesh of
$500$ points in the parametric domain $I$, while $E_R$ is defined as the relative error in $L^2(I)$ norm. The details on the considered BVP  are reported in the description of each numerical test.

\subsection{Exterior Dirichlet problem: parabola test}
\begin{figure}
\centering
\includegraphics[width=0.48\textwidth]{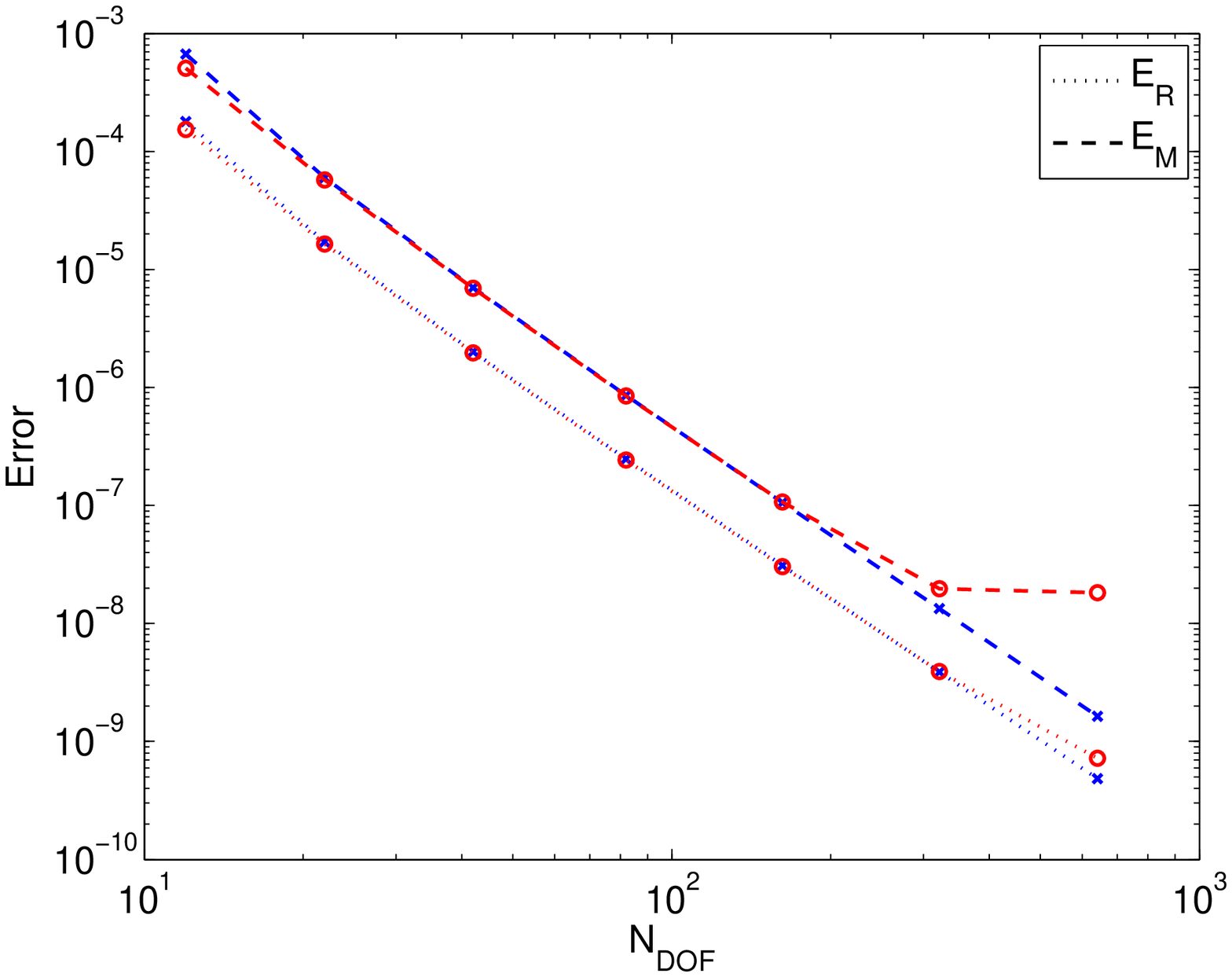} \includegraphics[width=0.48\textwidth]{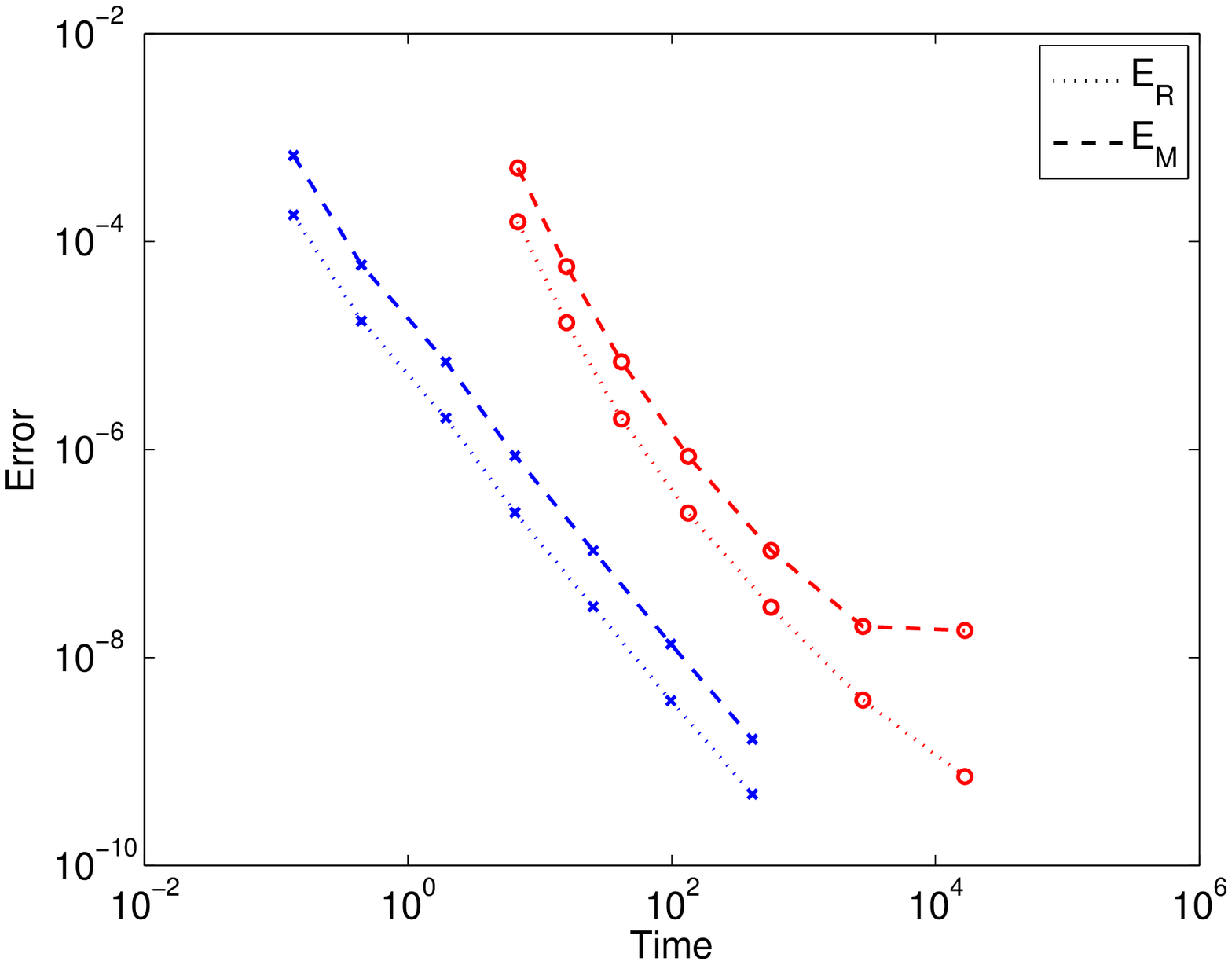}
\caption{Parabola test: Convergence of the errors, case $d=2$. On the left the plot of the errors with respect to the number of degrees of freedom $N_{DoF},$ on the right with respect to the assembly time. In blue with cross signs the new assembly, in red with circles the element-by-element strategy.  From the Figure on the left, the theoretical convergence of order $d+1=3$ can be noticed. From the Figure on the right, we can conclude that the new strategy is accurate as the previous and much more rapid. } \label{Parabola_conv}
\end{figure}

At first, let us consider the Dirichlet BVP defined in  (\ref{BVP_arco})
exterior to the arc of parabola $\Gamma=\{{\bf
x}=(x_1,x_2)\,|\,  x_1=t,\,x_2=1-t^2, \,t\in [-1,1]$\},
representable by means of quadratic B-splines related to the
extended knot vector
\begin{equation*}
 T_1= \left[ \begin{array}{c c c c c c}
 -1 &-1& -1& 1& 1& 1 \end{array} \right]
\end{equation*}

and to the control points $\Q_i,\, i=0,\,1,\,2$, whose
coordinates are collected in the following matrix:\\
\begin{equation*}
 Q=\left[\begin{array}{c c c}
-1 &0& 1\\0& 2& 0
\end{array}\right] \ .
\end{equation*}
Here, the Dirichlet datum $u_{\cal D}$ is given in such a way that the solution
of \eqref{prima_BIE} is explicitly known and
reads $\varphi(\x)=\sqrt{1+4 x_1^2}$.
The discrete counterpart of such BIE involves the computation of matrices $I_{K_1}$ and $I_{K_2}$ with, respectively, regular and weakly singular integrals, see equations \eqref{IK1e2}-\eqref{IK1e2_bis}. Also, the right-hand side is regular so that in our computations with the new assembly a B-spline weighted quadrature rule has been used.
The comparison reported in Table \ref{tabparall}, for different values of the parameter $h$, which uniformly decomposes the parameter interval
$[-1,1]$, involves $C^1$ quadratic B-spline basis used in the old (element-by-element, with $N_G=N_{prod}=32,\,N_{DE}=63$) and new ($N_{ref}= 1$ for both quadrature formulas)  implementations of IgA-SGBEM.\\
Together with DoF and spectral condition numbers of the associated matrices, we show the relative error $E_R$ in $L^2$ norm, the
absolute errors $E_M$ in maximum norm and the elapsed time for the IgA-SGBEM matrix generation. Note that the slightly different results concerning conditioning and accuracy are due to the different quadrature schemes adopted in the old and new assembly strategies. The elapsed time behaves in both cases as ${\cal O}(N_h^2)$ as expected, but the superiority of the new approach is self-evident. \\
Similar conclusions can be deduced looking at the results in the bottom part of Table \ref{tabparall} where the comparison has been done for increasing values of the
B-spline basis degree $d$ chosen for the approximation of the BIE solution. Here the elapsed time behaves in both cases as ${\cal O}(d^2)$.\\
 The same data are used in Figure \ref{Parabola_conv} where the convergence of errors are plotted with respect to the number of degrees of freedom and with respect to assembly time. In the first plot it can be seen that the solution obtained with the new strategy is accurate as the one calculated before, in the second plot it can be seen that the new strategy obtains very good errors using modest times. \\
  \begin{table} [thp]
\begin{center}
\resizebox{1\textwidth}{!}{
\begin{tabular} {|c c| c c c c | c c c c|}
\hline
\multicolumn{2}{|c|}{}& \multicolumn{4}{|c|}{element-by-element}
& \multicolumn{4}{c|}{new assembly} \\ \hline \hline
{$h$}& {$DoF$}& {$cond.$}& {$E_R$}&{$E_M$}&{time ($s$)}&{$cond.$}& {$E_R$}&{$E_M$}&{time ($s$)}\\
\hline
$1/5$ & $12$     & $7.19\cdot10^1$ & $1.54\cdot10^{-4}$     & $5.08\cdot10^{-4}$  & $6.80$& $8.56\cdot10^1$ & $1.79\cdot 10^{-4}$   & $6.67\cdot10^{-4}$   & $0.13$ \\
$1/10$ & $22$     & $1.87\cdot10^2$ & $1.65\cdot10^{-5}$     & $5.73\cdot10^{-5}$   & $15.96$ & $1.87\cdot10^2$ & $1.72\cdot 10^{-5}$   & $5.96\cdot10^{-5}$   & $0.44$ \\
$1/20$ & $42$     & $4.57\cdot10^2$ & $1.96\cdot10^{-6}$     & $6.92\cdot10^{-6}$ & $41.75$  & $4.50\cdot10^2$ & $2.01\cdot 10^{-6}$   & $6.98\cdot10^{-6}$   & $1.93$ \\
$1/40$ & $82$     & $1.01\cdot10^3$ & $2.43\cdot10^{-7}$     & $8.50\cdot10^{-7}$ & $134.29$ & $9.90\cdot10^2$ & $2.48\cdot 10^{-7}$   & $8.60\cdot10^{-7}$   & $6.47$ \\
$1/80$ & $162$   & $2.12\cdot10^3$ & $3.03\cdot10^{-8}$    & $1.07\cdot10^{-7}$  & $566.12$ & $2.07\cdot10^3$ & $3.08\cdot 10^{-8}$   & $1.06\cdot10^{-7}$   & $25.45$ \\
$1/160$ & $322$   & $4.34\cdot10^3$ & $3.90\cdot10^{-9}$    & $1.97\cdot10^{-8}$  & $2828.18$ & $4.25\cdot10^3$ & $3.85\cdot10^{-9}$    & $1.34\cdot10^{-8}$  & $98.87$ \\
$1/320$ & $642$   & $8.78\cdot10^3$ & $7.17\cdot10^{-10}$    & $1.82\cdot10^{-8}$  & $16659.45$ & $8.59\cdot10^3$ & $4.81\cdot10^{-10}$    & $1.64\cdot10^{-9}$  & $409.79$ \\
 \hline \hline
 {$d$}& {$DoF$}& {$cond.$}& {$E_R$}&{$E_M$}&{time ($s$)}&{$cond.$}& {$E_R$}&{$E_M$}&{time ($s$)}\\
 \hline
$2$ & $12$     & $7.19\cdot10^1$ & $1.54\cdot10^{-4}$     & $5.08\cdot10^{-4}$  & $6.80$& $8.56\cdot10^1$ & $1.79\cdot 10^{-4}$   & $6.67\cdot10^{-4}$   & $0.13$ \\
$3$ & $13 $     & $1.89\cdot10^2$ & $3.23\cdot 10^{-5}$   & $1.69\cdot10^{-4}$    & $11.93$ & $2.14\cdot10^2$ & $5.63\cdot 10^{-5}$   & $3.87\cdot10^{-4}$    & $0.21$\\
$4$ & $14$     & $5.09\cdot10^2$ & $1.80\cdot 10^{-5}$   & $8.01\cdot10^{-5}$    & $20.44$& $5.61\cdot10^2$ & $2.19\cdot 10^{-5}$   & $1.20\cdot10^{-4}$    & $0.25$\\
$5$ & $15 $     & $1.41\cdot10^3$ & $5.43\cdot 10^{-6}$   & $2.06\cdot10^{-5}$    & $34.27$ & $1.65\cdot10^3$ & $1.05\cdot 10^{-5}$   & $5.53\cdot10^{-5}$    & $0.31$\\
\hline
\end{tabular}}
\end{center}
\caption{Parabola test : comparison between the two different assembly strategy. On the top, for degree $d=2$ and various spacing $h$. On the bottom, for constant spacing $h=1/5$ and various degrees $d$.
On the first column we report calculated spectral condition number of the matrix. On the second and third we report the calculated errors: true solution is known, thus the relative error $E_R$ in $L^2$ norm and the absolute errors $E_M$ in maximum norm are computed. Finally the elapsed time in seconds for the IgA-SGBEM matrix generation is reported.\\
Notice that the error when $d$ increases decrees theoretically of a factor $h$, and this convergence is partially maintained from the methods, being the condition number increasing. By the other side the convergence in $h$ is as predicted in both cases, namely of order $d+1$. In all tested cases the condition number of the system is almost the same for the two assembly strategy and the overall time for assembly is strongly reduced by the new strategy. }
\label{tabparall}
\end{table}


\subsection{Interior Dirichlet problem: S-shaped closed domain test}
\begin{figure}
\centering
\includegraphics[width=0.48\textwidth]{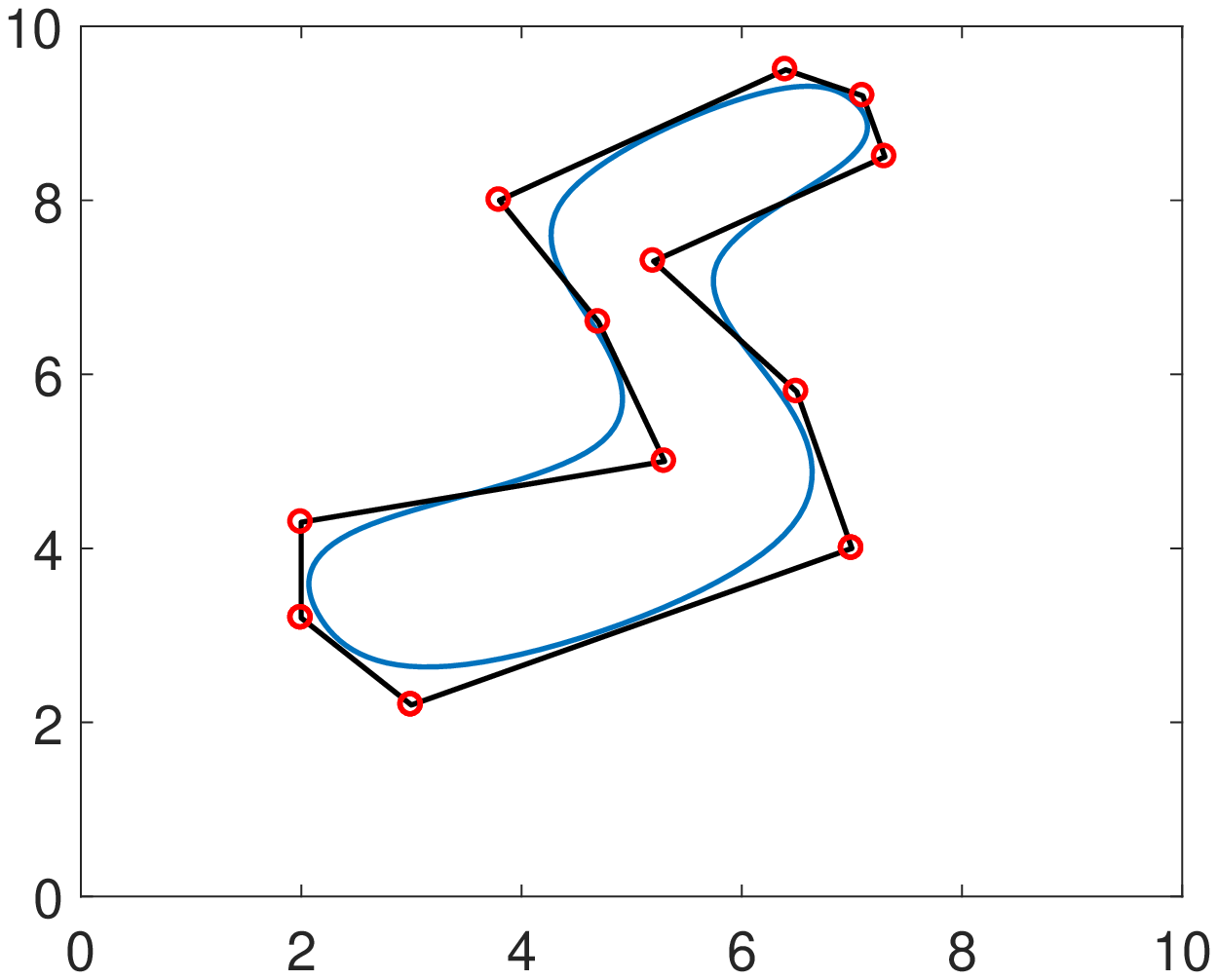} \includegraphics[width=0.48\textwidth]{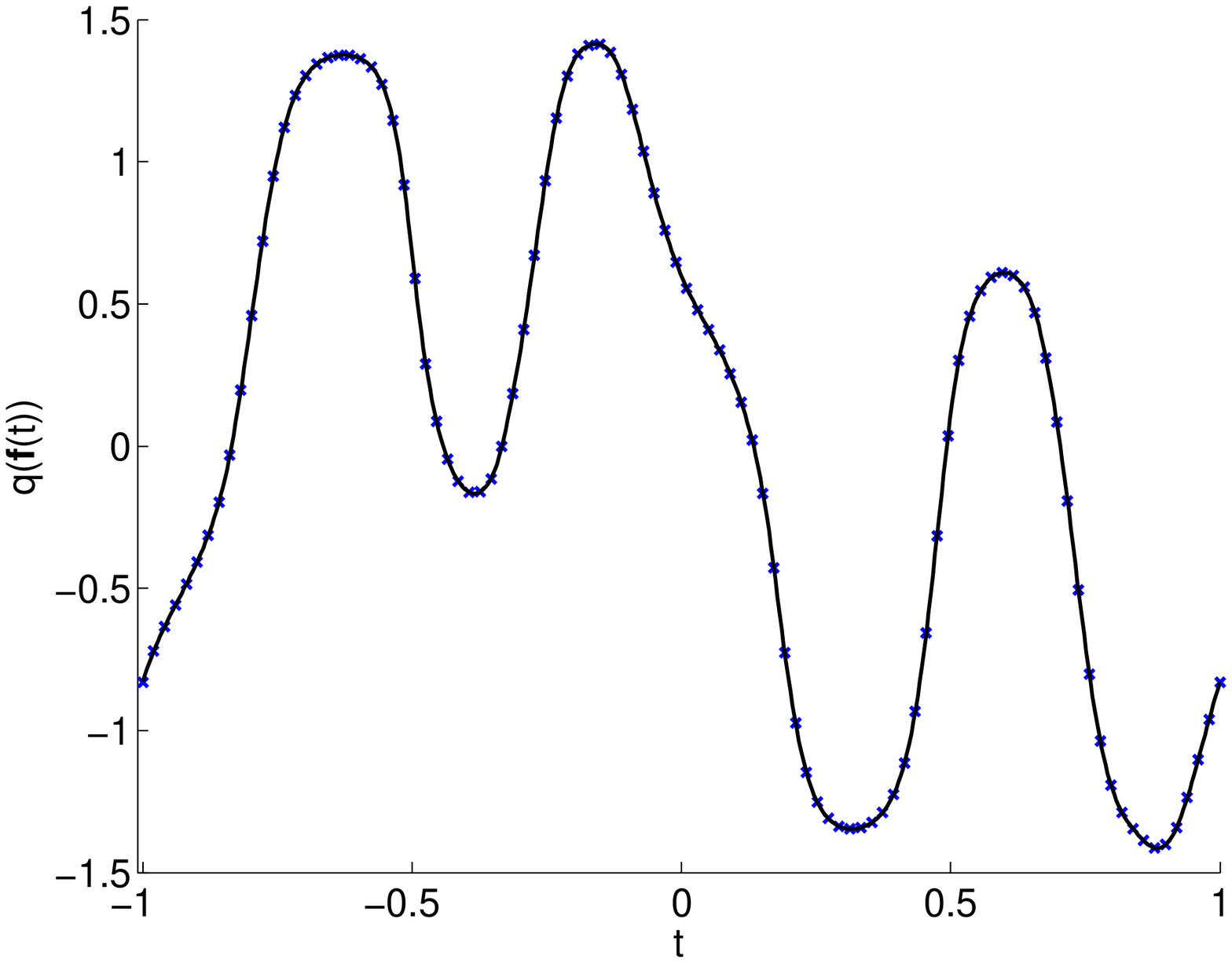}
\caption{S-shaped domain test: on the left the domain $\Omega$ with the control polygon
(black) used to represent its boundary in cubic B--spline form
with the  cyclic knot vector $T_2$ (in red circles the control points); on the right analytic solution (solid) and numerical one
obtained with $d=3, h=1/48$ (crosses).} \label{ESSE}
\end{figure}

We consider the interior Laplace problem (\ref{BVP_esse}) on the domain $\Omega$
shown in Figure \ref{ESSE}, equipped by Dirichlet boundary condition,
where $\Gamma:=\partial \Omega$ is described by cubic B-splines defined by the cyclic extended knot vector
$$T_2 \,=[\,-3/2\, :\, 1/6 \, :\, 3/2]\,,$$ and control points depicted in Figure \ref{ESSE}.
Since $u_{\cal D}=-(x_1+x_2)$ is chosen, the solution of \eqref{seconda_BIE} is explicitly known, it reads
$q(\x)=q(\f(t))=(f'_1(t)-f'_2(t))/\|\f'(t)\|_2$ and has $C^1([-1,1])$
regularity.\\
 Table \ref{tabesseall} compares results obtained
for this numerical test by the IgA-SGBEM element-by-element implementation ($N_G=N_{prod}=32,\,N_{DE}=63$) and the new one ($N_{ref}= 2$ for both the quadrature formulas).
 Also in this case, the new quadrature and assembly strategy reveals much faster than the old one.  Because of modest regularity of the solution the performed tests are limited to the case $d=3$.
The plot of
the approximate solution obtained from the new implementation with $h=1/48$
is depicted in
Figure \ref{ESSE}, where the analytical solution of \eqref{seconda_BIE} is also reported. \\
 This example is more challenging from the numerical point of view, due to modest regularity,
 the use of a direct approach and the oscillations of the solution, see Figure \ref{ESSE}.
 Moreover, we point out that the S-shaped domain, being not starred, is not good for some of the interior point-BEM methods introduced recently \cite{chen2016nurbs}.
 \begin{table} [htp]
\begin{center}
\resizebox{1\textwidth}{!}{
\begin{tabular} {|c c| c c c c | c c c c|}
\hline
\multicolumn{2}{|c|}{}& \multicolumn{4}{|c|}{element-by-element}
& \multicolumn{4}{c|}{new assembly} \\
 \hline \hline
 {$h$}& {$DoF$}& {$cond.$}& {$E_R$}&{$E_M$}&{time ($s$)}&{$cond.$}& {$E_R$}&{$E_M$}&{time ($s$)}\\
 \hline
$1/6$ & $15 $     & $3.39\cdot10^2$ & $1.12\cdot 10^{-1}$   & $1.62\cdot10^{-1}$    & $98.28$& $2.87\cdot10^2$ &   $1.14\cdot 10^{-1}$   & $2.33\cdot10^{-1}$   & $2.17$ \\
$1/12$ & $27$     & $8.20\cdot10^2$ &  $3.23\cdot 10^{-2}$   & $6.71\cdot10^{-2}$    & $351.41$& $9.19\cdot10^2$ &  $3.30\cdot 10^{-2}$   & $9.39\cdot10^{-2}$   & $5.97$ \\
$1/24$ & $51$     & $3.00\cdot10^3$ &  $3.98\cdot 10^{-3}$  & $1.33\cdot10^{-2}$   & $1240.44$ & $3.67\cdot10^3$ &  $4.03\cdot 10^{-3}$  & $1.97\cdot10^{-2}$   & $22.50$ \\
$1/48$ & $99$     & $8.75\cdot10^3$ &  $5.79\cdot 10^{-4}$   & $2.03\cdot10^{-3}$   & $4728.49$& $8.95\cdot10^3$ &  $6.32\cdot 10^{-4}$   & $3.11\cdot10^{-3}$   & $88.80$ \\
$1/96$ & $195$   & $2.09\cdot10^4$ &  $8.96\cdot 10^{-5}$   & $3.06\cdot10^{-4}$   & $18615.21$ & $2.15\cdot10^4$ &  $1.85\cdot 10^{-4}$   & $1.69\cdot10^{-3}$   & $372.18$ \\
\hline
\end{tabular}}
\end{center}
\caption{S-shaped domain test: comparison between the two different assembly strategies for degree $d=3$.  For the explanation of the columns see the previous table.
Notice that in all the tested cases the new assembly is accurate and efficient. }   \label{tabesseall}
\end{table}

\section{Conclusions and future work} \label{sect:concl}
In this paper we have presented a new strategy that gives fast assembly in IgA-SGBEM dealing with BVPs equipped by Dirichlet boundary conditions. The remaining analysis for the application to mixed problems will be covered in a forthcoming paper. \\
Note that the presented integration schemes can be easily used also in the context of collocation BEM, fixing the outer collocation node and using the new quadrature rules for the inner integration.\\
Moreover, the proposed method can be profitably used when treating non-linear transmission problems, in order to obtain the trace of the solution on the boundary, see \cite{Cal}.\\
An open issue consists in testing hierarchical spline spaces to develop an efficient adaptive isogeometric version of the scheme, see \cite{BG16}. The hierarchical approach could overcome the problem of dealing with non-regular curves.\\
The extension to 3D problems where the boundary surface can be described just by one open patch and the multi-patch case for the description of closed surfaces is currently under study.

\section*{Acknowledgements}
The support by Gruppo Nazionale per il  Calcolo Scientifico
(GNCS) of the Istituto Nazionale di Alta Matematica (INdAM)
through ``Progetti di ricerca'' program is gratefully
acknowledged. Giancarlo Sangalli
was partially supported by the European Research Council through the FP7 ERC Consolidator Grant No.616563 ''HIGEOM''.


\end{document}